\documentclass[10pt,reqno]{amsart}
\usepackage{amsfonts,amsthm,latexsym,amsmath,amssymb,amscd,amsmath}

\theoremstyle{plain}
   \newtheorem{theorem}{Theorem}[section]
   \newtheorem{proposition}[theorem]{Proposition}
   
   \newtheorem{corollary}[theorem]{Corollary}
   \newtheorem{conjecture}[theorem]{Conjecture}
   \newtheorem{problem}[theorem]{Problem}
   
\theoremstyle{definition}
   \newtheorem{definition}[theorem]{Definition}
   
\theoremstyle{remark}
\newtheorem{remark}[theorem]{Remark}
\newtheorem*{ack}{Acknowledgements}

\renewcommand{\Im}{\mathfrak{Im}}

\newcommand{\al}{\alpha}
\newcommand{\be}{\beta}
\newcommand{\la}{\lambda}
\newcommand {\bC} {\mathbb C}
\newcommand {\bR} {\mathbb R}

\newcommand {\ga} {\gamma}

\newcommand\pa {\partial}

\newcommand{\Si}{\mathfrak{S}}
\newcommand{\De}{\Delta}
\newcommand{\bN}{\mathbb N}
\newcommand{\HH}{\mathcal{H}}

\newcommand{\io}{\iota}
\newcommand{\RR}{\mathbb{R}}

\newcommand{\CC}{\mathbb{C}}

\newcommand{\cS}{\mathcal{S}}
\newcommand{\cT}{\mathcal{T}}
\newcommand{\cL}{\mathcal{L}}

\def\newop#1{\expandafter\def\csname #1\endcsname{\mathop{\rm
#1}\nolimits}}

\newop{diag}

\begin{document}
\numberwithin{equation}{section}
\title[Stable polynomials and mixed determinants]{Applications of 
stable polynomials to mixed determinants: Johnson's conjectures, unimodality 
and symmetrized Fischer products}
\author[J.~Borcea]{Julius Borcea}
\address{Department of Mathematics, Stockholm University, SE-106 91
Stockholm, Sweden}
\email{julius@math.su.se}
\author[P.~Br\"and\'en]{Petter Br\"and\'en}
\address{Department of Mathematics, Royal Institute of Technology, SE-100 44
Stockholm, Sweden}
\email{pbranden@math.kth.se}


\begin{abstract}
For $n \times n$ matrices  $A$ and $B$ define 
$$
\eta(A,B)
= \sum_{\cS} \det (A[\cS]) \det \!\left(B\!\left[\cS'\right]\right),   
$$
where the summation is over all subsets of $\{1,\ldots, n\}$, $\cS'$ is the 
complement of $\cS$, and $A[\cS]$ is the principal submatrix of $A$ with rows 
and columns indexed by $\cS$. We prove that if $A\ge 0$ and $B$ is Hermitian 
then (1) the polynomial $\eta(zA,-B)$ has all real roots (2) the latter 
polynomial has as many positive, negative and zero roots (counting 
multiplicities)
as suggested by the inertia of $B$ if $A>0$ and (3) for $1\le i\le n$ the 
roots of 
$\eta(zA\left[\{i\}'\right],-B\left[\{i\}'\right])$ interlace those of 
$\eta(zA,-B)$. Assertions (1)--(3) solve three important conjectures 
proposed by 
 C.~R.~Johnson 20 years ago. Moreover, we substantially extend these 
results to tuples of matrix pencils and real stable polynomials. 
In the process, we establish unimodality properties in the sense of 
majorization for the coefficients of 
homogeneous real stable polynomials and as an application, we derive 
similar properties for symmetrized Fischer products of positive definite 
matrices. We also obtain Laguerre type inequalities for characteristic 
polynomials of principal submatrices of arbitrary Hermitian matrices that
considerably generalize a certain subset of the 
Hadamard-Fischer-Koteljanskii inequalities for principal minors of
positive definite matrices. Finally, we propose Lax type problems 
for real stable polynomials and mixed determinants.
\end{abstract}

\maketitle


\section{Introduction and Main Results}\label{s1}

Let $A$ and $B$ be matrices of order $n$. For a  subset $\cS$ of $\{1,\ldots, n\}$ we denote by 
$A[\cS]$ the $|\cS| \times |\cS|$ principal submatrix of $A$ whose rows and columns are indexed by $\cS$. Define  
$$
\eta(A,B)= \sum_{\cS} \det (A[\cS] )\det\!\left(B\left[\cS'\right]\right),   
$$
where the summation is over all subsets of $\{1,\ldots, n\}$, $\cS'$ is the complement of $\cS$ and by convention $\det (A[\emptyset] )= \det (B[\emptyset])=1$. Note that 
$\eta(zI,-B)=\det(zI-B)$ is the characteristic polynomial of $B$. Thus, if $B$ is Hermitian then all zeros of $\eta(zI,-B)$ are real. Motivated by this 
fact and extensive numerical evidence, Johnson proposed the following conjecture around 1987.

\begin{conjecture}[\cite{J1}]\label{con1} 
If $A$ and $B$ are Hermitian matrices of the same order and $A$ is positive 
semidefinite then the polynomial 
$\eta(zA,-B)$ has all real zeros. 
\end{conjecture}

Let $\alpha_1 \leq \alpha_2 \leq \cdots \leq \alpha_n$ and 
$\beta_1 \leq \beta_2 \leq \cdots \leq \beta_m$ be the zeros of two 
univariate {\em hyperbolic} polynomials, i.e., two real polynomials with all
real zeros. These zeros are 
{\em interlaced} if they can 
be ordered so that either $\alpha_1 \leq \beta_1 \leq \alpha_2 \leq \beta_2 \leq 
\cdots$ or $\beta_1 \leq \alpha_1 \leq \beta_2 \leq \alpha_2 \leq 
\cdots$. By convention, we will also say that the zeros of the zero-polynomial interlace the zeros of any hyperbolic polynomial.

Conjecture~\ref{con1} was subsequently refined by Johnson \cite{J2} and 
also independently by Bapat \cite{B} to the following generalization of the classical Cauchy-Poincar\'e 
theorem stating that the eigenvalues of an $n \times n$ Hermitian matrix and those of any of its $(n-1)\times (n-1)$ principal submatrices interlace. 

\begin{conjecture}[\cite{J2}, \cite{B}]\label{con2}
If $A$ and $B$ are Hermitian $n\times n$ matrices and $A$ is positive semidefinite then the zeros of 
$\eta(zA[\{j\}'],-B[\{j\}'])$ interlace those of $\eta(zA,-B)$ for  $1 \leq j \leq n$ (provided that $\eta(zA,-B)$ is not identically zero).
 \end{conjecture}

We note that as originally formulated in \cite{J1,J2} the assumptions of 
Conjectures \ref{con1}--\ref{con2} actually 
required $A$ to be positive 
definite. It is however easy to see that if true, these conjectures should
in fact hold for any positive semidefinite $n\times n$ matrix $A$ and 
Hermitian $n\times n$ matrix $B$. In \cite{J2} it was 
further conjectured that if $A$ is positive definite and $B$ is
Hermitian then $\eta(zA,-B)$ has as many positive, negative and zero roots
as the inertia of $B$ would suggest. More precisely, given a univariate 
hyperbolic polynomial
$p(z)$ denote by $\io_{+}(p(z))$, $\io_{-}(p(z))$
and $\io_{0}(p(z))$ the number of positive, negative and zero roots of $p(z)$,
respectively. For any such polynomial we may then define its ``inertia''
as $\io(p(z))=(\io_{+}(p(z)),\io_{0}(p(z)),\io_{-}(p(z)))$.

\begin{conjecture}[\cite{J2}]\label{con3}
If $A$ and $B$ are Hermitian $n\times n$ matrices and $A$ is positive definite then $\io(\eta(zA,-B))=\io(\eta(zI,-B))\,(=\io(\det(zI-B)))$.
\end{conjecture}

A successful approach to Conjectures \ref{con1}--\ref{con2} has so far remained
elusive and only some special cases have been dealt with in the literature. For instance, Conjecture~\ref{con1} has been verified for $n=3$
in \cite{Ru} (albeit in a rather complicated way) and in \cite{B}, where 
 Conjectures \ref{con1} and \ref{con2} were shown to hold when both $A$ 
and $B$ are tridiagonal. These results were subsequently generalized in
\cite{F1,F2} to matrices whose graph is a tree or a cycle. Except for 
considerable empirical evidence (cf.~\cite{J2}) not much seems to be known
concerning Conjecture \ref{con3}. 

In \S \ref{s2} we prove Conjectures \ref{con1}--\ref{con3} in full 
generality and we significantly extend Conjectures~\ref{con1}--\ref{con2} 
to mixed determinants for tuples of matrix pencils and multivariate real 
stable polynomials (Definition~\ref{def1}, Definition~\ref{def-md} and Theorem~\ref{master}). Mixed  determinants are similar to -- but not the same as -- mixed discriminants and mixed volumes that originate from the classical work of Minkowski, Weyl, and later, Alexandrov (\cite{Gu2}, 
\cite[Chap.~2.3]{Hor}; see also Remark~\ref{re-m-disc}).  

In \S \ref{s3} we establish log-concavity and monotonicity properties
in the sense of majorization for the coefficients of 
real stable homogeneous polynomials.
As a consequence, we obtain similar properties for symmetrized
Fischer products associated with positive definite matrices 
(Theorem~\ref{t-log} and 
Corollaries~\ref{cor1}--\ref{c-gen}). 

In \S \ref{new-s3} we apply the theory of real stable polynomials
developed in \cite{BBS1,BBS2,Br} and the present paper 
to derive Laguerre type inequalities for 
characteristic polynomials of principal submatrices of arbitrary Hermitian 
matrices (Theorem~\ref{t-kot} and Corollary~\ref{c-kot}). These polynomial
inequalities vastly generalize a certain subset of the classical 
Hadamard-Fischer-Koteljanskii inequalities for principal minors of
positive definite matrices. It is interesting to note that this same 
subset of determinantal inequalities actually appears in many other contexts,
such as the study of the 
``Principal Minor Assignment Problem'' for real symmetric matrices 
\cite{HS} (see Remark~\ref{r-hol-sturm}).

Finally, in \S \ref{s5} we propose and discuss Lax type 
representations for real stable (homogeneous) polynomials by means of
mixed determinants (Problems~\ref{pb1}--\ref{pb2}). In view of the close 
connections between real stable polynomials and G\aa rding hyperbolic 
polynomials (cf.~Proposition~\ref{pro-g} in \S \ref{s5}) these problems 
are natural steps -- albeit from a somewhat different perspective than 
Helton-Vinnikov's \cite{HS} -- toward finding higher dimensional analogs of 
the Lax conjecture (now a theorem, see 
Theorems~\ref{dub}--\ref{laxlike} in \S \ref{s5}). 
In this context we should also mention that in recent work \cite{Gu1,Gu2}
Gurvits has successfully used G\aa rding hyperbolic polynomials to generalize
and reprove in a unified manner a number of classical conjectures in matrix
theory and real algebraic geometry, including the van der Waerden and 
Schrijver-Valiant conjectures as well as Alexandrov-Fenchel type inequalities
for mixed discriminants and mixed volumes.

\begin{ack}
We would like to thank the three anonymous referees for their useful comments.
It is a pleasure to thank the American Institute of Mathematics for 
hosting the ``P\'olya-Schur-Lax Workshop'', May-June 2007,  
on topics pertaining
to this project, see \cite{BBCV}.
\end{ack}

\section{Proofs of Conjectures~\ref{con1}--\ref{con3}}\label{s2}

As we shall see, Conjectures \ref{con1}--\ref{con3} follow from a general theory developed 
in an ongoing series of papers \cite{BBS1,BBS2, Br}, where we study generalizations of the notion of "real-rootedness" to several variables. 

Recall from \S \ref{s1} that a nonzero univariate polynomial with real coefficients is said to be 
{\em hyperbolic} if all its zeros are real. A univariate polynomial 
$f(z)$ with complex coefficients is called {\em stable} if 
$f(z) \neq 0$ for all $z \in \CC$ with $\Im(z) >0$. Hence a univariate 
polynomial with real coefficients is stable if and only if it is hyperbolic. 

These classical concepts have several natural extensions to multivariate 
polynomials. Here as well as in \cite{BBS1} we are concerned with the most 
general notion of this type, which may be defined as follows.

\begin{definition}\label{def1}
A polynomial $f\in\bC[z_1,\ldots,z_n]$ is called 
{\em stable} 
if  $f(z_1, \ldots, z_n) \neq 0$ for all $n$-tuples $(z_1, \ldots, z_n) \in \CC^n$ with $\Im(z_j) >0$, $1\leq j \leq n$.  If in addition $f$ has real coefficients it will be referred  to as {\em real stable}. We denote by 
$\HH_n(\CC)$ and $\HH_n(\RR)$ the set of stable and real stable polynomials, respectively, and by $\HH_n^+(\RR)$ the subset of $\HH_n(\RR)$ consisting of
polynomials with all nonnegative coefficients.
\end{definition}

Thus 
$f$ is stable (respectively, real stable) if and only if for all $\alpha \in \RR^n$ and $v \in \RR_+^n$ -- where we let as usual $\RR_+=(0,\infty)$ -- the 
univariate polynomial $f(\alpha + vt)$ is stable (respectively, hyperbolic).  

Here are a few examples of (real) stable polynomials from various areas:
\begin{enumerate} 
\item Let $G=(V,E)$ be a finite connected graph. Then the {\em spanning tree polynomial} is the polynomial 
$$
T(G;z)= \sum_{T}z^T, \quad z^T:=\prod_{e \in T}z_e, 
$$ 
where the sum is over all spanning trees $T \subseteq E$ and $z_e$, $e \in E$, are (commuting) variables. This polynomial is real stable, see \cite{Br,COSW}.
\item Let $A$ be an $n \times n$ Hermitian matrix and let $Z=\diag(z_1, \ldots, z_n)$ be a diagonal matrix of (commuting) variables. Then $\det(A+Z)$ is a real stable polynomial, see Proposition~\ref{pencil} below. 
\item If $z\mapsto\sum_{k=0}^n a_k z^k$ is a stable or hyperbolic univariate
polynomial then 
$$
(z_1,\ldots,z_n)\mapsto \sum_{k=0}^n a_k \binom{n}{k}^{-1}e_k(z_1, \ldots, z_n)
$$
is a (multivariate) polynomial in $\HH_n(\CC)$, where $e_k(z_1, \ldots, z_n)$, $0\leq k \leq n$, are the elementary symmetric functions in $n$ variables. This follows from the Grace-Walsh-Szeg\"o coincidence theorem, see, e.g., \cite[Theorem 3.4.1b]{RS}. 
\end{enumerate}

The following are two basic properties of stable polynomials that we need to prove Conjectures \ref{con1}--\ref{con3}. 

\begin{proposition}\label{misc}
Let $f \in \HH_n(\CC)$. Then 
\begin{enumerate} 
\item $ {\partial f}/{\partial z_j} \in \HH_n(\CC)\cup\{0\}$ for $1 \leq j \leq n$, 
\item $f(z_1, \ldots, z_{j-1}, \alpha_j, z_{j+1}, \ldots, z_n) \in \HH_{n-1}(\CC)\cup\{0\}$ for $1\leq j \leq n$ and any $\alpha_j \in \CC$ with $\Im(\alpha_j) \geq 0$. 
\end{enumerate}
\end{proposition}

The first property follows from the Gauss-Lucas theorem, see also \cite{BBS1} where the authors characterize all finite order linear differential operators with polynomial coefficients that preserve stability. The second property is deduced by applying the multivariate version of Hurwitz' theorem stated below with
$D= \{z \in \CC : \Im(z)>0\}^n$ and $f_k(z_1,\ldots,z_n)= f(z_1, \ldots, z_{j-1}, \alpha_j+z_j/k, z_{j+1}, \ldots, z_n)$, $k\in\bN$.

\begin{theorem}[Hurwitz' theorem]\label{mult-hur}
Let $D$ be a domain (open connected set) in $\CC^n$ and suppose that $\{f_k\}_{k=1}^\infty$ is a sequence of nonvanishing 
analytic functions on $D$ that converge to $f$ uniformly on compact subsets of $D$. Then $f$ is either 
nonvanishing on $D$ or else identically zero. 
\end{theorem}

For a proof of Theorem~\ref{mult-hur} we refer to 
\cite[Footnote 3, p.~96]{COSW}. 
A key ingredient  in our proofs of Conjectures \ref{con1}--\ref{con3} is the following simple albeit fundamental proposition, see \cite[Proposition 2]{BBS1}. 

\begin{proposition} \label{pencil}
 Let $A_j$, $1\le j\le n$, be positive semidefinite $m \times m$ matrices 
and  let $B$ be a Hermitian  $m \times m$ matrix. Then the polynomial 
\begin{equation}\label{determ}
f(z_1,\ldots, z_n)=\det\!\left(\sum_{j=1}^nz_jA_j+ B\right)
\end{equation}
is either real stable or identically zero.   
\end{proposition}

\begin{proof}
By a standard continuity argument using Hurwitz' theorem it suffices to prove 
the result only in the case when  all matrices $A_j$, $1\le j\le n$, are 
positive definite. Set $z(t)=\alpha+\la t$ with
$\alpha=(\al_1,\ldots,\al_n)\in \RR^n$, 
$\la=(\la_1,\ldots,\la_n)\in \RR_+^n$ and $t\in \bC$. Note that 
$P:=\sum_{j=1}^n\la_jA_j$
is positive definite and thus it is invertible and has a square root (recall
that $\bR_+=(0,\infty)$, cf.~the paragraph following Definition~\ref{def1}). 
Then 
\begin{equation*}
f(z(t))=\det(P)\det(tI+P^{-1/2}HP^{-1/2}),
\end{equation*}
where $H:=B+\sum_{j=1}^n\al_jA_j$ is a Hermitian $m\times m$ matrix. 
Therefore $f(z(t))$ is a polynomial in $t$ which is a 
constant multiple of the characteristic polynomial of a
Hermitian matrix and so it must have all real zeros.
\end{proof}

\begin{definition}\label{def-md} 
The {\em mixed determinant} of an $m$-tuple $(A_1,\ldots,A_m)$ of $n\times n$ matrices is given by 
\begin{equation}\label{def-eta}
\eta(A_1,\ldots,A_m)=\sum_{(\cS_1,\ldots,\cS_m)}\det(A_1[\cS_1])\cdots \det(A_m[\cS_m]),
\end{equation}
where the summation is taken over all ordered partitions of $\{1, \ldots, n\}$ into $m$ parts i.e., all $m$-tuples $(\cS_1,\ldots,\cS_m)$ of 
pairwise disjoint subsets of 
$\{1,\ldots,n\}$ satisfying $\cS_1\cup\cdots\cup\cS_m=\{1,\ldots,n\}$.
\end{definition} 

We begin by proving a considerable generalization of Conjectures \ref{con1} 
and \ref{con2}.

\begin{theorem}\label{master}
Let $\ell, m,n \geq 1$ be integers. For $1 \leq j \leq m$ let
$$
\cL_j := \cL_j(z_1, \ldots, z_\ell) = \sum_{k=1}^\ell{A_{jk}z_k} + B_j, 
$$
where $A_{jk}$, $1 \leq k \leq \ell$, are positive semidefinite $n \times n$ 
matrices and $B_j$ is a Hermitian $n \times n$ matrix. Then 
$$
\eta(\cL_1,\ldots,\cL_m) \in \HH_\ell(\RR)\cup\{0\}.
$$
Moreover, if $v$ is a new variable then
$$
\eta(\cL_1,\ldots,\cL_m)+v\eta\left(\cL_1\left[\{j\}'\right],\ldots,\cL_m\left[\{j\}'\right]\right) \in \HH_{\ell+1}(\RR)\cup\{0\},  \quad 1\leq j \leq n.
$$
\end{theorem}

\begin{proof} 
 Let $W=\diag(w_1,\ldots,w_n)$, where $w_1, \ldots, w_n$ are variables, and 
let 
$$g_j(z_1,\ldots,z_\ell, w_1, \ldots, w_n):= \det(W+\cL_j),\quad 1\le j\le m.$$
By Proposition~\ref{pencil} we have that 
 $g_j \in \HH_{\ell+n}(\RR)$ for $1\le j\le m$. Clearly, the negative of the inverse of a complex number is in the open upper half-plane if and only if the number itself is in the open upper half-plane. Hence for $1\le i\le n$ we have
that $\Im(w_i)>0 \Leftrightarrow -\Im\!\left(w_i^{-1}\right)>0$ and therefore  
 $$
 \det(I-W\cL_j)=(-1)^nw_1\cdots w_n g_j(z_1, \ldots, z_\ell, -w_1^{-1}, \ldots, -w_n^{-1}) \in \HH_{\ell+n}(\RR)
 $$
for $1\le j\le m$. Letting $w=(w_1,\ldots,w_n)$ and $w^\cS=\prod_{i\in\cS}w_i$
we thus get 
\begin{equation}\label{prod}
\begin{split}
f(z_1,\ldots,z_\ell,w_1,\ldots,w_n):&=\prod_{j=1}^m\det(I-W\cL_j)\\
&= \prod_{j=1}^m\Big(\sum_{\cS}\det (\cL_j[\cS]) (-1)^{|\cS|}w^{\cS}\Big) \in \HH_{\ell+n}(\RR)
\end{split}
\end{equation}
since products of stable polynomials are stable. When expanding both sides of identity \eqref{prod} in powers of 
$w=(w_1,\ldots,w_n)$ we see that the coefficient of $w_1w_2\cdots w_n$ in $f(z_1,\ldots,z_\ell,w_1,\ldots,w_n)$ is $(-1)^n\eta(\cL_1,\ldots,\cL_m)$, i.e., 
\begin{equation}\label{mixed-tayl}
\eta(\cL_1,\ldots, \cL_m)= (-1)^n \frac {\partial^nf}{ \partial w_1 \cdots \partial w_n} \Big|_{w_1=\cdots=w_n=0}. 
\end{equation}
By Proposition~\ref{misc} the latter polynomial belongs to $\HH_{\ell}(\RR) \cup \{0\}$, which proves the first part of the theorem.
%

Let $V_j$ be the $n \times n$ diagonal matrix with all entries equal to zero but the $j$-th diagonal entry which is $v$. By a straightforward computation we get 
$$
\eta(V_j, \cL_1,\ldots,\cL_m) =  \eta(\cL_1,\ldots,\cL_m)+v\eta(\cL_1[\{j\}']\ldots,\cL_m[\{j\}']) 
$$
and by the above this is either a real stable polynomial in the variables $v,z_1,\ldots,z_\ell$ or identically zero, which settles the second part of the theorem.
\end{proof}

\begin{remark}\label{re-m-disc}
Consider an $n$-tuple $(A_1,\ldots,A_n)$ of $n\times n$ complex matrices. Then $\det(t_1A_1+\ldots+t_nA_n)$ is a homogeneous polynomial of degree $n$ in $t_1,\ldots,t_n$. The number (sometimes used with a normalizing factor $1/n!$) 
\begin{equation}\label{eq-mx-disc}
\mathcal{M}(A_1,\ldots,A_n):=\frac{\partial^n}{\partial t_1\ldots\partial t_n}\det(t_1A_1+\ldots+t_nA_n)
\end{equation}
is called the {\em mixed discriminant} of $A_1,\ldots,A_n$, see, e.g., \cite{Gu1,Gu2}. By Definition~\ref{def-md} and~\eqref{mixed-tayl} we know that the mixed determinant $\eta(A_1,\ldots,A_n)$ may also be expressed as a Taylor coefficient of a certain polynomial. Note though that despite the similarity between~\eqref{mixed-tayl} and~\eqref{eq-mx-disc} the notions of mixed determinant and mixed discriminant are actually quite different.
\end{remark}

We shall also use the classical Hermite-Biehler theorem \cite[Theorem 6.3.4]{RS}. 

\begin{theorem}[Hermite-Biehler theorem]
Let $h=f+ig \in \CC[z] \setminus \{0\}$, where $f,g \in \RR[z]$. Then $h$ is stable if and only if $f$ and $g$ are hyperbolic polynomials whose zeros 
interlace and $f'(z)g(z)-f(z)g'(z)\ge 0$ for all $z\in \RR$.
\end{theorem}

\begin{proof}[Proof of Conjectures \ref{con1} and \ref{con2}]
Conjecture \ref{con1} follows from Theorem~\ref{master} by letting $m=2$,  
$\cL_1=zA$ and $\cL_2=-B$. From Theorem~\ref{master} we also get that 
$$
\eta(zA,-B)+v\eta(zA[\{j\}'],-B[\{j\}']) \in \HH_2(\RR) \cup\{0\}.
$$
Setting $v=i$ yields 
$$
\eta(zA,-B)+i\eta(zA[\{j\}'],-B[\{j\}']) \in \HH_1(\CC) \cup\{0\}, 
$$
which by the Hermite-Biehler theorem proves Conjecture~\ref{con2}.
\end{proof}

\begin{proof}[Proof of Conjecture \ref{con3}]
Let $A,B$ be Hermitian $n\times n$ matrices such that $A>0$, i.e., $A$ is positive definite. We begin by proving that $\io_0(\eta(zA,-B))$ is equal to the nullity $\nu:=\dim(\ker(B))$ of $B$.  Clearly,  
$\nu=\min\{ |\cS|: \det(B[\cS']) \neq 0 \}$. 

It is evident from the definition of $\eta(zA,-B)$ that 
$\io_0(\eta(zA,-B))\geq \nu
$,  
so we just have to prove the reverse inequality. If $\nu=0$ then 
$\det(B) \neq 0$ and since the constant term of $\eta(zA,-B)$ is $(-1)^n\det(B)$ we also have 
$\io_0(\eta(zA,-B))=0$. If $\nu>0$ let $\cS=\{s_1,\ldots, s_\nu\}$ be such that $\det(B[\cS'])\neq 0$. Define a sequence  of polynomials by    
$p_0(z)= \eta(zA,-B)$ and $p_i(z)= \eta(zA[\{s_1,\ldots, s_i\}'],-B[\{s_1,\ldots,s_i\}'])$ for 
$1 \leq i \leq \nu$. Since Conjectures~\ref{con1}--\ref{con2} are valid we know that 
the zeros of $p_{i-1}(z)$ and $p_{i}(z)$ interlace for $1 \leq i \leq \nu$. In particular, we have that 
$$
\io_0(\eta(zA,-B))=\io_0(p_0) \leq \io_0(p_1)+1 \leq \io_0(p_2)+2 \leq \cdots \leq \io_0(p_\nu)+\nu.  
$$
By the choice of $\cS$ we know that $p_\nu(0)=(-1)^{n-\nu}\det(B[\cS'])\neq 0$ so that 
$\io_0(p_\nu)=0$ and thus $\io_0(\eta(zA,-B)\leq \nu$, as was to be proved. 

By the above, we have 
$\io_{0}(\eta(zI,-B))=\io_{0}(\eta(zA,-B))=\nu$. 
To complete the proof, note first that 
for all Hermitian matrices $B$ and positive definite matrices $A$ one has 
$\deg(\eta(zA,-B))=n$ since 
$\det(A)>0$. Suppose now that 
$\io_{+}(\eta(zA,-B))\neq \io_{+}(\eta(zI,-B))$ for some $A>0$. Let 
$A(t) = (1-t)I+tA$, $t\in [0,1]$, be a 
homotopy between $I$ and $A$. It is clear that $A(t)>0$, $t\in [0,1]$. 
By the above, we know that $\io_{0}(\eta(zA(t),-B))=\nu$ for 
$t\in [0,1]$. So 
we may write $\eta(zA(t),-B)= z^\nu p(t,z)$,
where $p(t,z)$ is a hyperbolic polynomial of degree $n-\nu$ and $p(t,0)\neq 0$ for all $t\in [0,1]$.
Now $\io_{+}(p(0,z))\neq \io_{+}(p(1,z))$ and so at least one zero of $p(t,z)$ 
has to pass through the origin as $t$ runs from 
$0$ to $1$. By Hurwitz' theorem we must therefore have $p(T,0)=0$ for some 
$T\in [0,1]$,  contrary 
to the assumption that $p(t,0) \neq 0$ whenever $0\leq t \leq 1$.
\end{proof}  

\section{Unimodality Properties for Stable Homogeneous Polynomials \\ and
Symmetrized Fischer Products}\label{s3}

Given a matrix $A$ of order $d$ and $0\le k\le d$ define the $k$-th symmetrized
Fischer product associated with $A$ by
$$
S_k(A)=\sum_{|\cS|=k}\det(A[\cS]) \det \!\left(A\!\left[\cS'\right]\right)
$$
and the corresponding $k$-th average Fischer product
$$\overline{S}_k(A)=\binom{d}{k}^{-1}S_k(A).$$

As a first application of Theorem~\ref{master} we establish the following 
result. We note that part (a) below was initially conjectured in 
\cite{J2} and later proved in \cite{BJ} by means of immanantal inequalities
derived from the theory of generalized matrix functions. 

\begin{corollary}\label{cor1}
If $A$ is a positive semidefinite $d\times d$ matrix then its average Fischer products satisfy 
\begin{itemize}
\item[(a)] $\overline{S}_k(A)\le \overline{S}_\ell(A)$ for $0\le k\le \ell \le \left[\frac{d}{2}\right]$, 
\item[(b)] $\overline{S}_k(A)^2 \geq \overline{S}_{k-1}(A)\overline{S}_{k+1}(A)$ for $1\le k\le d-1$, 
\item[(c)] $ \frac {\overline{S}_1(A)}{\det(A)} \geq \sqrt{\frac {\overline{S}_2(A)}{\det(A)}} \geq \sqrt[3]{\frac {\overline{S}_3(A)}{\det(A)}} \geq \cdots \geq \sqrt[d]{\frac {\overline{S}_d(A)}{\det(A)}}=1$ whenever $A>0$.
\end{itemize}
\end{corollary}

\begin{proof}
By Hurwitz' theorem we may assume that $A>0$. 
Note that
$$\eta(zA,-A)=\sum_{k=0}^{d}(-1)^{d-k}S_k(A)z^k.$$

From the above proofs of  Conjectures
\ref{con1} and \ref{con3} we have that all roots of the polynomial 
$\eta(zA,-A)$ are real and positive. Then by applying Newton's inequalities, 
see, e.g., \cite{ineq}, we immediately get part (b), which combined with 
Maclaurin's inequalities \cite{ineq} yields part (c). 

Since $A$ is positive definite we have $\overline{S}_k(A)>0$ for  $0\leq k \leq d$. 
It is well known that property (b)  and positivity imply that the sequence  $\left\{\overline{S}_k(A)\right\}_{k=0}^d$ is unimodal, i.e., it is weakly increasing until it reaches a peak after which it  is weakly decreasing (see, e.g., \cite[p. 137]{wilf}). Since 
$\left\{\overline{S}_k(A)\right\}_{k=0}^d$ is also
symmetric -- that is, $\overline{S}_k(A)=\overline{S}_{d-k}(A)$ for $0\leq k \leq d$ -- part (a) follows. 
\end{proof}

The monotonicity and log-concavity properties given by  
Corollary \ref{cor1} (a)-(b) may actually be viewed as special cases of 
a more general phenomenon, namely unimodality for the coefficients of 
real stable homogeneous polynomials that we proceed to describe. Let 
$\{e_i\}_{i=1}^n$ denote the standard basis in $\bR^n$ and recall the 
{\em majorization preordering} $\prec$ 
on $n$-tuples of real numbers (cf.~\cite{ineq,MO}): if 
$x=(x_1,\ldots,x_n)\in\bR^n$ and $y=(y_1,\ldots,y_n)\in\bR^n$ one says that
$y$ majorizes $x$ and writes $x\prec y$ provided that
$$\sum_{i=1}^kx_{[i]}\le \sum_{i=1}^{k}y_{[i]},\quad 1\le k\le n,$$
with equality for $k=n$, where for $w=(w_1,\ldots,w_n)\in\bR^n$ we let 
$(w_{[1]}\ge \ldots\ge w_{[n]})$ denote its decreasing rearrangement. 
An $n$-tuple $x=(x_1,\ldots,x_n)$ is said to be a {\em transfer} (or 
{\em pinch}) of another $n$-tuple $y=(y_1,\ldots,y_n)$ if there exist
$i\in\{1,\ldots,n-1\}$ and 
$t\in \left(0,\frac{y_{[i]}-y_{[i+1]}}{2}\right]$ such that $x_{[i]}=y_{[i]}-t$,
$x_{[i+1]}=y_{[i+1]}+t$ and $x_{[k]}=y_{[k]}$ for $k\in\{1,\ldots,n\}\setminus
\{i,i+1\}$. According to a 
well-known theorem of Hardy-Littlewood-P\'olya and Muirhead, if $x,y\in\bR^n$
then $x\prec y$ if and only if $x$ may be obtained from $y$ by a finite 
number of pinches (cf.~\cite{ineq,MO}).

For $\al=(\al_1,\ldots,\al_n)\in\bN^n$ with $\al_1+\cdots+\al_n=d$ we define as usual
$$\binom{d}{\al}=\frac{d!}{\al_1!\cdots\al_n!}.$$
Given a real homogeneous polynomial of degree $d$ 
\begin{equation}\label{poly}
f(z)=\sum_{|\al|=d}a(\al)z^{\al}\in\bR[z_1,\ldots,z_n]
\end{equation}
(where we use the standard multi-index notation $z=(z_1,\ldots,z_n)$, 
$|\al|=\sum_{i=1}^n\al_i$ and $z^{\al}=\prod_{i=1}^{n}z_i^{\al_i}$) we 
normalize its coefficients by setting
$$\hat{a}(\al)=\binom{d}{\al}^{-1}a(\al).$$

\begin{theorem}\label{t-log}
Let $f$ be a real homogeneous polynomial of degree $d$ 
as in~\eqref{poly}. If $f\in\HH_n(\bR)$ and 
$\al=(\al_1,\ldots,\al_n)\in\bN^n$ is such that $|\al|=d$ and 
$\al_i\ge \al_j>0$ then
$$\hat{a}(\al+k(e_i-e_j))^2\ge 
\hat{a}(\al+(k-1)(e_i-e_j))\hat{a}(\al+(k+1)(e_i-e_j))$$
for $-\al_i+1\le k\le \al_j-1$.
\end{theorem}

\begin{proof}
If $d\le 1$ there is nothing to prove so we may assume that $d\ge 2$. 
Now in order to establish the inequality stated in the theorem it is clearly 
enough 
to show that if $\ga=(\ga_1,\ldots,\ga_n)\in\bN^n$ is such that $|\ga|=d$
and $\ga_i\ga_j>0$ then
\begin{equation}\label{eq-extra}
\hat{a}(\ga)^2\ge \hat{a}(\ga-e_i+e_j)\hat{a}(\ga+e_i-e_j).
\end{equation}
Let $\be\in\bN^n$ with $|\be|=d-2$, set $\ga=\be+e_i+e_j$ and define a polynomial $g$ in two variables $z_i$ and $z_j$ by
$$g(z_i,z_j)=\frac{\pa^{\be}f}{\pa z^{\be}}\bigg|_{z_k=0,\,k\neq i,j}.$$
An elementary computation yields
$$g(z_i,z_j):=\frac{d!}{2}\left[\hat{a}(\be+2e_i)z_i^2+
2\hat{a}(\be+e_i+e_j)z_iz_j+\hat{a}(\be+2e_j)z_j^2\right].$$
Since $f\in\HH_n(\bR)$ it follows from the properties of (real) stable polynomials given
in Proposition~\ref{misc} that $g\in\HH_2(\bR)\cup\{0\}$ and thus 
$g(z_i,1)\in\HH_1(\bR)\cup \{0\}$. By evaluating the 
discriminant of the latter univariate (hyperbolic) polynomial we get
$$\hat{a}(\be+e_i+e_j)^2\ge \hat{a}(\be+2e_i)\hat{a}(\be+2e_j),$$
which is the same as~\eqref{eq-extra}.
\end{proof}

\begin{corollary}\label{c-monot}
Let $f$ be a real homogeneous polynomial of degree $d$ 
as in~\eqref{poly} which we further assume to be symmetric in all its variables and with at least
one positive coefficient. 
If $f\in\HH_n(\bR)$ and 
$\al=(\al_1,\ldots,\al_n)\in\bN^n$ is such that $|\al|=d$ and 
$\al_i\ge \al_j$ then
$$\hat{a}(\al+k(e_i-e_j))\le \hat{a}(\al+\ell(e_i-e_j))$$
for $-\al_i\le k\le \ell\le \left[\frac{\al_j-\al_i}{2}\right]$ and
$$\hat{a}(\al+k(e_i-e_j))\ge \hat{a}(\al+\ell(e_i-e_j))$$
for $\left[\frac{\al_j-\al_i}{2}\right]\le k\le \ell\le \al_j$.
\end{corollary}

\begin{proof}
Since $f$ is a real stable polynomial which is also homogeneous, all its 
nonzero coefficients must have the same phase (that is, the quotient of any two
such coefficients is a positive number), see
\cite[Theorem 6.1]{COSW}. By assumption, $f$ has at least one positive
coefficient and therefore all its (nonzero) coefficients are
positive, i.e., $f\in\HH_n^+(\bR)$. It then follows from the log-concavity 
property established in Theorem~\ref{t-log} that the sequence
$$\left\{\hat{a}(\al+k(e_i-e_j))\right\}_{k=-\al_i}^{\al_j}$$
is unimodal (cf.~the proof of Corollary~\ref{cor1}). Moreover, this sequence
also satisfies 
$$\hat{a}(\al+k(e_i-e_j))
=\hat{a}(\al+(\al_i+\al_j-k)(e_i-e_j))$$ 
for $-\al_i\le k\le \al_j$ since 
$f$ is assumed to be symmetric in the variables $z_1,\ldots,z_n$. This 
proves the desired result. 
\end{proof} 

\begin{corollary}\label{c-major}
Assume that the real homogeneous polynomial $f$ given 
by~\eqref{poly} is symmetric in all its variables and has at least one positive coefficient. 
If $f\in\HH_n(\bR)$ and $\al,\be\in\bN^n$ are such that 
$\be\prec \al$
(in the sense of majorization) then $\hat{a}(\al)\le \hat{a}(\be)$.
\end{corollary}
 
\begin{proof}
As in the proof of Corollary~\ref{c-monot} we deduce from 
\cite[Theorem 6.1]{COSW} and from 
the fact that $f$ is a real stable homogeneous 
polynomial that
$f\in\HH_n^+(\bR)$. Now by the aforementioned standard properties of the 
majorization preordering
it is enough to check the assertion only for $\al,\be\in\bN^n$ such that
$\be$ is a pinch of $\al$. Since the coordinates of 
$\al=(\al_1,\ldots,\al_n)$ and those of 
any of its pinches should be nonnegative integers this is in turn equivalent 
to showing that if $i,j$ are such that $\al_i> \al_j$ then 
$\hat{a}(\al)\le \hat{a}(\al-e_i+e_j)$. The latter inequality is an immediate
consequence of Corollary~\ref{c-monot}.
\end{proof}

Corollary \ref{cor1} above deals with symmetrized Fischer products for
positive definite $n\times n$ matrices 
corresponding to partitions of $n$ with only two parts. 
Theorem~\ref{t-log} and Corollaries~\ref{c-monot}--\ref{c-major} allow us to 
extend 
Corollary \ref{cor1} to symmetrized Fischer products corresponding to 
arbitrary partitions of $d$. These products are defined as follows: let 
$\al=(\al_1,\ldots,\al_n)\in\bN^n$ be such that $|\al|=d$. Given a matrix 
$A$ of order $d$ set
$$S_{\al}(A)=\sum_{(\cS_1,\ldots,\cS_n)}
\prod_{i=1}^{n}\det \!\left(A\!\left[\cS_i\right]\right),$$
where the summation is taken over all ordered partitions 
 $(\cS_1,\ldots,\cS_n)$ of $\{1,\ldots,d \}$ into $n$ parts such that
$|\cS_i|=\al_i$, $1\le i\le n$. The corresponding average Fischer
product is then given by 
$$\overline{S}_{\al}(A)=\binom{d}{\al}^{-1}S_{\al}(A).$$
Note that both $S_{\al}(A)$ and $\overline{S}_{\al}(A)$ have a natural 
$\Si_n$-invariance property, where $\Si_n$ denotes the symmetric group on 
$n$ elements.
Indeed, if $\al=(\al_1,\ldots,\al_n)\in\bN^n$ and 
$\pi(\al)=(\al_{\pi(1)},\ldots,\al_{\pi(n)})$ then 
$$S_{\pi(\al)}(A)=S_{\al}(A)\text{ and }
\overline{S}_{\pi(\al)}(A)=\overline{S}_{\al}(A)\text{ for all }\pi\in\Si_n.$$

\begin{corollary}\label{c-gen}
If $A$ is a positive semidefinite $d\times d$ matrix then
\begin{itemize}
\item[(a)] $\overline{S}_{\al}(A)\le \overline{S}_{\be}(A)$ for any
$\al,\be\in\bN^n$ such that $|\al|=|\be|=d$ and $\be\prec\al$,
\item[(b)] 
whenever $\al=(\al_1,\ldots,\al_n)\in\bN^n$
is such that $|\al|=d$, $\al_i\ge \al_j>0$ and $-\al_i+1\le k\le \al_j-1$
one has
$$\overline{S}_{\al+k(e_i-e_j)}(A)^2\ge 
\overline{S}_{\al+(k-1)(e_i-e_j)}(A)\,\overline{S}_{\al+(k+1)(e_i-e_j)}(A).$$
\end{itemize}
\end{corollary}

\begin{proof}
In view of Hurwitz's theorem we may assume that $A$ is positive definite. 
From Definition~\ref{def-md} we then deduce that
$$f(z_1,\ldots,z_n):=\eta(z_1A,\ldots,z_nA)\in\bR[z_1,\ldots,z_n]$$
is a symmetric homogeneous polynomial of degree $d$ with all positive 
coefficients. By Theorem~\ref{master} one has $f\in\HH_n(\bR)$. Applying
Theorem~\ref{t-log} and Corollary~\ref{c-major} to $f$ one immediately
gets parts (b) and (a), respectively. 
\end{proof}

\begin{remark}
Corollary \ref{c-gen} (a) was proved using other methods in 
\cite[Theorem~1]{BJ} (see the paragraph preceding Corollary \ref{cor1} 
above).
\end{remark} 

\section{Laguerre Type Extensions to Hermitian Matrices of the 
Hadamard-Fischer-Koteljanskii Inequalities}\label{new-s3}

Any positive definite $n\times n$ matrix $A=(a_{ij})$ satisfies a number of 
well-known determinantal inequalities \cite{FJ,HJ} including classical ones
such as
\begin{align*}
&\text{Hadamard: }\quad &&\det (A)\le \prod_{i=1}^{n}a_{ii},\\
&\text{Fischer: }\quad &&\det(A)\le \det(A[\cS])\det(A[\cS']),\\
&\text{Koteljanskii: }\quad &&\det(A[\cS\cup\cT])\det(A[\cS\cap\cT])
\le \det(A[\cS])\det(A[\cT]),\\
\intertext{where $\cS,\cT\subseteq \{1,\ldots,n\}$. These inequalities are 
actually valid even for other classes of matrices, such as totally positive
matrices 
and $M$-matrices (cf.~{\em op.~cit.}). However, as noted
in e.g.~\cite{FJ} Hadamard's, Fischer's and Koteljanskii's inequalities are 
essentially equivalent for the class of positive definite matrices. We will 
now show that a certain subset of the Hadamard-Fischer-Koteljanskii
inequalities holds for arbitrary Hermitian matrices. In fact we establish a 
more general result involving characteristic polynomials of principal 
submatrices of Hermitian matrices which is reminiscent of yet another
classical inequality, namely}
&\text{Laguerre: }\quad &&f(z)f''(z)\le f'(z)^2,\quad z\in\bR,
\end{align*}
whenever $f\in\HH_1(\bR)$ (see, e.g., \cite[Lemma 5.4.4 and p.~179]{RS}). More precisely, we prove
the following result.

\begin{theorem}\label{t-kot}
Let $A$ be a Hermitian $n\times n$ matrix and $\cS,\cT\subseteq \{1,\ldots,n\}$
be such that $|\cS\cap\cT|=|\cS|-1=|\cT|-1$. For any $z\in\bR$ one has
\begin{multline*}
\det\!\left(zI_{|\cS\cup\cT|}-A[\cS\cup\cT]\right)
\det\!\left(zI_{|\cS\cap\cT|}-A[\cS\cap\cT]\right)\\
\le\det\!\left(zI_{|\cS|}-A[\cS]\right)\det\!\left(zI_{|\cT|}-A[\cT]\right).
\end{multline*}
\end{theorem}
Setting $z=0$ in Theorem~\ref{t-kot} and taking into account that the sign on both sides will be plus due to the sizes of the principal submatrices we get:
 
\begin{corollary}\label{c-kot}
If $A$ is a Hermitian $n\times n$ matrix and $\cS,\cT\subseteq \{1,\ldots,n\}$
are such that $|\cS\cap\cT|=|\cS|-1=|\cT|-1$ then
$$\det(A[\cS\cup\cT])\det(A[\cS\cap\cT])
\le \det(A[\cS])\det(A[\cT]).$$
\end{corollary}

\begin{remark}\label{r-hol-sturm}
It is interesting to note that the condition imposed on the cardinalities
of the index sets $\cS,\cT$ in Theorem~\ref{t-kot} and Corollary~\ref{c-kot}
is identical with condition (12) in Theorem 5 of \cite{HS} dealing with 
the ``Principal Minor Assignment Problem'' for certain types of vectors in 
$\bR^{2^n}$ and real symmetric $n\times n$ matrices.
\end{remark}

Theorem~\ref{t-kot} is an easy consequence of the characterization
of real stable polynomials that was recently obtained in \cite[Theorem 10]{Br}
(see also \cite[Theorem 27]{BBS1}). In the case of multi-affine polynomials
(i.e., of degree at most one in each variable) the criterion
for real stability established in \cite{Br} may be formulated as 
follows:

\begin{theorem}\label{rs-crit}
A multi-affine polynomial $f\in\bR[z_1,\ldots,z_n]$ is real stable if and only
if $\De_{ij}(f)(x_1,\ldots,x_n)\ge 0$ for all $(x_1,\ldots, x_n)\in\bR^n$ and 
$1\le i,j\le n$,
where 
$$\De_{ij}(f)=\frac{\partial f}{\partial z_i}\cdot\frac{\partial f}
{\partial z_j}-\frac{\partial^2 f}{\partial z_i\partial z_j}\cdot f
=-f^2\frac{\partial^2 }{\partial z_i\partial z_j}\left[\log |f|\right].$$
\end{theorem}

\begin{proof}[Proof of Theorem~\ref{t-kot}]
Since principal submatrices of Hermitian matrices are themselves 
Hermitian we may assume 
without loss of generality that $\cS\cup\cT=\{1,\ldots,n\}$.
Let then $i,j\in\{1,\ldots,n\}$ be such that 
$\cS\cup\{i\}=\cT\cup\{j\}=\{1,\ldots,n\}$
and consider the polynomial
$$f(z_1,\ldots,z_n)=\det(\text{diag}(z_1,\ldots,z_n)-A).$$
Proposition~\ref{pencil} implies that $f\in\HH_n(\bR)$ and thus 
$\De_{ij}(f)(z,\ldots,z)\ge 0$ for any $z\in \bR$ by Theorem~\ref{rs-crit}. 
A straightforward computation now yields
\begin{multline*}
\De_{ij}(f)(z,\ldots,z)=\det\!\left(zI_{|\cS|}-A[\cS]\right)
\det\!\left(zI_{|\cT|}-A[\cT]\right)\\
-\det\!\left(zI_{|\cS\cup\cT|}-A[\cS\cup\cT]\right)
\det\!\left(zI_{|\cS\cap\cT|}-A[\cS\cap\cT]\right),
\end{multline*}
which proves the theorem.
\end{proof}

A further interesting consequence of Theorem~\ref{t-kot} is the following.

\begin{corollary}\label{c-inert}
If $A$ is a Hermitian $n\times n$ matrix with $\det(A)=0$ then 
$$\det\!\left(A\left[\{i\}'\right]\right)
\det\!\left(A\left[\{j\}'\right]\right)\ge 0$$ 
for all $i,j\in\{1,\ldots,n\}$.
\end{corollary}

\section{Lax Type Problems for Stable Polynomials}\label{s5}

Recall that a   
homogeneous polynomial $f\in\bR[z_1,\ldots, z_n]$  is said to be {\em (G\aa rding) hyperbolic} with respect to a given 
vector $e \in \RR^n$ if $f(e) \neq 0$ and for all vectors $\al \in \bR^n$  the univariate polynomial $f(\alpha + et)\in\bR[t]$ has all real zeros. 
As is well known, such polynomials play an 
important role in e.g.~the theory of partial differential operators. 
It turns out that real stable polynomials and multivariate homogeneous 
hyperbolic polynomials are closely related objects, as noted in 
\cite[Proposition 1]{BBS1}:

\begin{proposition}\label{pro-g}
Let $f(z_1,\ldots, z_n)$ be a polynomial of degree $d$ with real coefficients and let 
$p(z_0,z_1, \ldots, z_n)$ be the (unique) homogeneous polynomial of degree $d$ such that 
$p(1,z_1, \ldots, z_n)=f(z_1,\ldots,z_n)$. Then $f$ is real stable if and only if 
$p$ is hyperbolic with respect to every vector $e\in\bR^{n+1}$ such that $e_0=0$ and $e_i>0$, $1 \leq i \leq n$. 
\end{proposition}

In 1958 Lax conjectured that any (G\aa rding) hyperbolic polynomial in three variables 
admits a certain determinantal representation \cite{Lax}. 
The Lax conjecture has recently been verified in 
\cite{LPR}:

\begin{theorem}[\cite{LPR}]\label{dub}
A homogeneous polynomial $f \in \RR[x,y,z]$  is hyperbolic of degree $d$ with respect to the 
vector $e=(1,0,0)$ if and only if there exist real symmetric $d\times d$ 
matrices $B,C$ such that $f(x,y,z)=f(e)\det(xI+yB+zC)$ for any $x,y,z\in\bR$. 
\end{theorem}

\begin{remark}
The proof of Theorem~\ref{dub} given in \cite{LPR} essentially follows from
the results of \cite{HV}. It is worth mentioning that a preliminary 
result along these lines was earlier obtained in \cite[Theorem 6.4]{Du}. 
\end{remark}

Using Theorem~\ref{dub} a converse to 
Proposition~\ref{pencil} in the case $n=2$ was established in 
\cite[Theorem 11]{BBS1}. More precisely, the following natural analog of the 
Lax conjecture for real stable polynomials in two variables was obtained
in {\em loc.~cit.}  

\begin{theorem}[\cite{BBS1}]\label{laxlike}
Any real stable polynomial in two variables $x, y$ can be written as 
$\pm\det(xA+yB+C)$ where $A$ and $B$ are positive semidefinite matrices 
and $C$ is a symmetric matrix of the same order. 
\end{theorem}   

It is however well known that the analog of the Lax conjecture fails in the case of four or 
more variables. The following modified version of the higher 
dimensional Lax conjecture has recently been proposed in \cite{HV}
(see also \cite[Conjecture 1]{BBS1}):

\begin{conjecture}\label{HV}
Let $P(x_0,x_1,\ldots,x_m)$
be a real homogeneous polynomial hyperbolic with respect to
$c=(c_0,c_1,\ldots,c_m) \in {\mathbb R}^{m+1}$ and $L$ be a real linear form 
in 
$x_0,x_1,\ldots,x_m$ with
$L(c) \neq 0$. Then there exists an integer $N$ such that
$$
L(x_0,x_1,\ldots,x_m)^N P(x_0,x_1,\ldots,x_m) 
= \det(x_0 A_0 + x_1 A_1 + \ldots +x_m A_m)
$$
for some real symmetric matrices $A_0,A_1,\ldots,A_m$
with $c_0 A_0 + c_1 A_1 + \ldots + c_m A_m > 0$. 
\end{conjecture}

A natural question in this context is whether any real stable (homogeneous) 
polynomial admits a Lax type determinantal representation. To formulate
precise versions of this question 
let $\ell, m,n \geq 1$ be integers. For $1 \leq j \leq m$ define
matrix pencils 
\begin{equation}\label{mt-pen}
\cL_j := \cL_j(z_1, \ldots, z_\ell) = \sum_{k=1}^\ell{A_{jk}z_k} + B_j, 
\end{equation}
where $A_{jk}$, $1 \leq k \leq \ell$, are positive semidefinite $n \times n$ 
matrices and $B_j$ is a Hermitian $n \times n$ matrix. Then by 
Theorem~\ref{master} we know that 
$\eta(\cL_1,\ldots,\cL_m) \in \HH_\ell(\RR)\cup\{0\}$.

\begin{problem}\label{pb1}
Is the converse of (the first part of) Theorem~\ref{master} true, namely:
if $f$ is a real stable polynomial of degree $n$ in $\ell$ variables 
then there exist a positive integer $m$ and matrix
pencils $\cL_j$, $1\le j\le m$, of the 
form~\eqref{mt-pen} such that $f=\eta(\cL_1,\ldots,\cL_m)$?
\end{problem}

Note that by Theorem~\ref{laxlike} the answer to Problem~\ref{pb1} is 
affirmative (at least) in the case $\ell=2$. The homogeneous version of
Problem~\ref{pb1} is as follows.

\begin{problem}\label{pb2}
Let $f$ be a real stable homogeneous polynomial of degree $n$ in 
$\ell$ variables. Is it true that there exist a 
positive integer $m$ and matrix
pencils $\cL_j$, $1\le j\le m$, of the 
form~\eqref{mt-pen} with $B_j=0$, $1\le j\le m$, 
such that $f=\eta(\cL_1,\ldots,\cL_m)$?
\end{problem}

\end{document}